# ASYMPTOTIC BEHAVIOR OF EDGE-REINFORCED RANDOM WALKS

By Franz Merkl and Silke W. W. Rolles

*University of Munich and Technische Universität München*

In this article, we study linearly edge-reinforced random walk on general multi-level ladders for large initial edge weights. For infinite ladders, we show that the process can be represented as a random walk in a random environment, given by random weights on the edges. The edge weights decay exponentially in space. The process converges to a stationary process. We provide asymptotic bounds for the range of the random walker up to a given time, showing that it localizes much more than an ordinary random walker. The random environment is described in terms of an infinite-volume Gibbs measure.

**1. Introduction.** *Edge-reinforced random walk* on a locally finite undirected graph is the following process: Every edge is assigned a weight which changes with time. Initially, all weights equal a constant $a$. The random walker starts at a vertex **0**. At every time, the random walker jumps to a neighboring vertex with probability proportional to the weight of the traversed edge at that time. Each time an edge is traversed, its weight is increased by 1.

This model was introduced by Diaconis in [1] and [2]. The process is partially exchangeable. Already in 1980, Diaconis and Freedman [3] proved for the more general class of partially exchangeable processes a representation as a mixture of Markov chains, provided the process is recurrent.

In the late 1980s, Diaconis asked whether edge-reinforced random walk on $\mathbb{Z}^d$ is recurrent. Except for $d = 1$, this problem is still unsolved. On an infinite binary tree, Pemantle [8] showed a phase transition in the recurrence and transience behavior of edge-reinforced random walk. For general *finite graphs*, Coppersmith and Diaconis [1] found an explicit description for the limiting fraction of time spent at the edges. Their result was extended by Keane and Rolles [6]. In [9], Rolles showed that a class of models can be represented as a mixture of *reversible* Markov chains. This result applies in particular to edge-reinforced random walk on any finite graph. Edge-reinforced random walks were used in [4] to provide natural Bayesian priors for reversible Markov chains.

In [7], we proved that the edge-reinforced random walk on the ladder $\mathbb{Z} \times \{1,2\}$ is recurrent for all initial edge weights $a > 3/4$. This result was generalized by one of the authors in [10] to graphs of the type $\mathbb{Z} \times G$ and $\mathbb{N}_0 \times G$, where $G$ is a finite tree, provided that the initial weights are sufficiently large.

In this article, we examine the asymptotic behavior of these edge-reinforced random walks on the infinite ladder $\mathbb{N}_0 \times G$ in much more detail beyond recurrence.

*Formal description of the model and notation.* We consider edge-reinforced random walk on a graph $\overline{G} = (\overline{V}, \overline{E})$. The vertex set $\overline{V}$ is of the type $\overline{V} = \mathbb{N}_0 \times V$ with a finite set $V$, $|V| \geq 2$. The set $V$ is assumed to be the vertex set of a finite tree $G = (V, E)$. Two vertices $(i, v), (i', v') \in \overline{V}$ are connected by an edge in $\overline{E}$ iff $i = i'$ and $v, v'$ are connected by an edge in $E$, or $|i - i'| = 1$ and $v = v'$. The edges of $\overline{G}$ are undirected. The edge-reinforced random walker starts at level 0 of $\overline{G}$, that is, in a vertex $\mathbf{0} = (0, v)$.

Furthermore, we assume the initial weights $a$ to be sufficiently large. More quantitatively, we assume $a > a_{\min}$, where $a_{\min} = 3/4$ if $V = \{1, 2\}$, and otherwise $a_{\min} = a_{\min}(G)$ denotes the lower bound specified in formula (1.7) of [10]. Optimizing the lower bound for $a$ is not treated in this paper.

The edge-reinforced random walk on $\overline{G}$ is formally defined as follows: Let $X_t : \overline{V}^{\mathbb{N}_0} \to \overline{V}$ denote the canonical projection on the $t$th coordinate; $X_t$ is interpreted as the location of the random walker at time $t$. For $t \in \mathbb{N}_0$, we define $w_t(e) : \overline{V}^{\mathbb{N}_0} \to \mathbb{R}_+$, the weight of edge $e$ at time $t$, recursively as follows:

$$(1.1) \qquad w_0(e) := a \qquad \text{for all } e \in \overline{E},$$

$$(1.2) \qquad w_{t+1}(e) := \begin{cases} w_t(e) + 1, & \text{for } e = \{X_t, X_{t+1}\} \in \overline{E}, \\ w_t(e), & \text{for } e \in \overline{E} \setminus \{\{X_t, X_{t+1}\}\}. \end{cases}$$

The distribution $P_\mathbf{0}$ of the edge-reinforced random walk is a probability measure on $\overline{V}^{\mathbb{N}_0}$, defined by

$$(1.3) \qquad X_0 = \mathbf{0}, \qquad P_\mathbf{0}\text{-a.s.},$$

$$(1.4)\ P_\mathbf{0}[X_{t+1} = v | X_i, i = 0, 1, \ldots, t] = \begin{cases} \dfrac{w_t(\{X_t, v\})}{\sum_{\{e \in \overline{E} : X_t \in e\}} w_t(e)}, & \text{if } \{X_t, v\} \in \overline{E}, \\ 0, & \text{otherwise.} \end{cases}$$

For a vertex $(i, v) \in \overline{V}$, we set $|(i, v)| := |i|$, and we abbreviate $v_i := (i, v)$. If $e = \{u, v\}$ is an edge in the finite graph $G$, we set $e_i := \{u_i, v_i\}$. For an edge $e = \{u, v\} \in \overline{E}$, we denote by $|e|$ its distance from level 0: we set $|e| := \min\{|u|, |v|\}$. Constants are denoted by $c_i$, $i \geq 1$. They keep their meaning throughout the whole article.

## 2. Results.

2.1. *Position of the random walker at large times.* Typically, at time $t$, simple random walk is located at a distance of the order $\sqrt{t}$ from its starting point. In contrast to this fact, the typical location of edge-reinforced random walk at time $t$ does not go to infinity as $t$ grows. In fact, the location of the reinforced random walk at time $t$ is stochastically bounded by a random variable with exponential tails, uniformly for all times $t$. This is the claim of the following theorem:

THEOREM 2.1 (Uniform exponential tails for the location of the random walk). *There exist constants $c_1, c_2 > 0$, depending only on $\overline{G}$ and $a$, such that for all $t, n \in \mathbb{N}_0$, the following bound holds:*

$$(2.1) \qquad P_{\mathbf{0}}(|X_t| \geq n) \leq c_1 e^{-c_2 n}.$$

As a consequence, up to time $t$, the edge-reinforced random walk can travel at most a distance of the order $\ln t$:

COROLLARY 2.2 (Range of the random walk path up to a given time). *There exists a constant $c_3 = c_3(\overline{G}, a) > 0$ such that $P_{\mathbf{0}}$-a.s.,*

$$(2.2) \qquad \max_{s=0,\ldots,t} |X_s| \leq c_3 \ln t \qquad \text{for all } t \text{ large enough.}$$

Simple random walk does not converge to an equilibrium distribution. Reinforcement makes this behavior change drastically. As is shown in the following theorem, the law of the location of the reinforced random walk tends to an equilibrium distribution as time grows.

However, the graph $\overline{G}$ has a chessboard structure. Half of the vertices can only be reached in an even number of steps, and the other half only in an odd number of steps. Therefore, we can only expect limit theorems for the reinforced random walk restricted to all even or all odd times. For these two restricted processes, we have indeed the following limit theorem:

THEOREM 2.3 (Convergence to equilibrium). *As $t \to \infty$, the distributions of $X_{2t}$ and $X_{2t+1}$ converge in the following sense: There exist probability functions $\mu_{\text{even}}$ and $\mu_{\text{odd}}$ on the vertex set $\overline{V}$, such that for all vertices $v \in \overline{V}$, the following hold:*

$$(2.3) \qquad \lim_{t \to \infty} P_{\mathbf{0}}(X_{2t} = v) = \mu_{\text{even}}(v),$$

$$(2.4) \qquad \lim_{t \to \infty} P_{\mathbf{0}}(X_{2t+1} = v) = \mu_{\text{odd}}(v).$$

*With the constants $c_1, c_2 > 0$ from Theorem* 2.1, *one has for all vertices* $v \in \overline{V}$:

$$\mu_{\text{even}}(v) \le c_1 e^{-c_2|v|} \quad \text{and} \quad \mu_{\text{odd}}(v) \le c_1 e^{-c_2|v|}. \tag{2.5}$$

Let $\overline{V}_{\text{even}}$ and $\overline{V}_{\text{odd}}$ denote the set of vertices $v \in \overline{V}$ which can be reached in an even and odd number of steps, respectively, by the random walker starting from $\mathbf{0}$. The measures $\mu_{\text{even}}$ and $\mu_{\text{odd}}$ are supported on $\overline{V}_{\text{even}}$ and $\overline{V}_{\text{odd}}$, respectively.

2.2. *Random walk in random environment representation.* Edge-reinforced random walk on any *finite* graph can be represented as a unique mixture of reversible Markov chains. This is shown in Theorem 3.1 of [9].

Recall that the transition probabilities of any irreducible reversible Markov chain on any graph $(\mathcal{V}, \mathcal{E})$ can be described by weights $x = (x_e)_{e \in \mathcal{E}}$, $x_e \ge 0$, on the edges of the graph; the probability to traverse an edge is proportional to its weight. More precisely, denoting the distribution of the Markov chain induced by the edge weights $x$ with starting vertex $v$ by $Q_{v,x}$, one has

$$Q_{v,x}(X_{t+1} = u' | X_t = u) = \frac{x_{\{u,u'\}}}{x_u}, \tag{2.6}$$

where here and in the following we set

$$x_u := \sum_{e \in \mathcal{E} : u \in e} x_e. \tag{2.7}$$

Edge-reinforced random walk on the *infinite* graph $\overline{G}$ can also be represented as a unique mixture of reversible Markov chains. Moreover, the corresponding random weights $(x_e)_{e \in \overline{E}}$ are summable. This is shown by the following theorem.

We introduce the infinite simplex $\Delta := \{(x_e)_{e \in \overline{E}} \in (0,1)^{\overline{E}} : \sum_{e \in \overline{E}} x_e = 1\}$.

THEOREM 2.4 (Mixture of positive recurrent Markov chains). *The edge-reinforced random walk on the* infinite *ladder $\overline{G}$ can be represented as a unique mixture of the reversible Markov chains $Q_{v,x}$. Even more, the mixing measure is supported on positive recurrent Markov chains. Hence, there is a unique probability measure $\mathbb{Q}$ on $\Delta$ such that*

$$P_{\mathbf{0}}(A) = \int_\Delta Q_{\mathbf{0},x}(A)\,\mathbb{Q}(dx) \tag{2.8}$$

*is valid for all events $A \subseteq \overline{V}^{\mathbb{N}_0}$.*

Let $\overline{G}^{(n)} = (\overline{V}^{(n)}, \overline{E}^{(n)})$ be the restriction of $\overline{G}$ to the finite vertex set $\{0, 1, \ldots, n\} \times V$. On $\overline{G}^{(n)}$, we can also describe the mixing measure as a measure on the space of weights. Let $x^{(n)} = (x_e^{(n)})_{e \in \overline{E}^{(n)}}$ denote the edge weights, normalized such that $\sum_{e \in \overline{E}^{(n)}} x_e^{(n)} = 1$. We denote the mixing measure for the weights $x^{(n)}$ by $\mathbb{Q}^{(n)}$.

Let $e_0^*$ be a fixed edge in the rung at level 0, incident to the starting vertex **0**, that is, $e_0^* = \{(0, u), (0, v)\}$ for some $\{u, v\} \in E$ and $\mathbf{0} \in e_0^*$. The following theorem shows that the edge weights decay exponentially in space with probabilities exponentially close to 1, even when we normalize the weights by dividing by $x_{e_0^*}^{(n)}$. This is true on the infinite ladder, but also on any finite ladder, uniformly in the size of the ladder.

THEOREM 2.5 (Exponential decay of the edge weights). *There exist positive constants $c_4, c_5, c_6$ depending only on $\overline{G}$ and $a$ such that for all $n \in \mathbb{N}$ and all edges $e$ of $\overline{G}^{(n)}$, we have*

$$(2.9) \qquad \mathbb{Q}^{(n)}(x_e^{(n)} > x_{e_0^*}^{(n)} e^{-c_4|e|}) \leq c_5 e^{-c_6|e|},$$

*uniformly in $n$. On the infinite ladder $\overline{G}$, we have the similar bound*

$$(2.10) \qquad \mathbb{Q}(x_e > x_{e_0^*} e^{-c_4|e|}) \leq c_5 e^{-c_6|e|}$$

*for all $e \in \overline{E}$.*

As a corollary, we obtain that the weights decay exponentially in space almost surely, even uniformly from a certain (random) point on. We can use any fixed edge $f$ to normalize the weights.

COROLLARY 2.6 (Exponential decay of the edge weights). *Let $f \in \overline{E}$ be a fixed edge. There exists a positive constant $c_4(\overline{G}, a) > 0$ with the following property: For $\mathbb{Q}$-almost all $x \in \Delta$, there exists a (random) $n \in \mathbb{N}_0$, such that for all edges $e$ with $|e| \geq n$, one has*

$$(2.11) \qquad x_e \leq e^{-c_4|e|/2} x_f.$$

Let $\pi = (v_0 = \mathbf{0}, v_1, \ldots, v_k)$ be a finite path in $\overline{G}$. Conditioned on $(X_0, \ldots, X_k) = \pi$, the shifted process $(X_{k+t})_{t \in \mathbb{N}_0}$ is again an edge-reinforced random walk starting at $v_k$ with initial edge weights given by

$$(2.12) \qquad a + \sum_{i=1}^{k} 1_{\{\{v_{i-1}, v_i\} = e\}}, \qquad e \in \overline{E}.$$

The shifted process can also be represented as a unique mixture of reversible Markov chains:

THEOREM 2.7 (ERRW conditioned on a finite path). *For any finite path $\pi = (v_0 = \mathbf{0}, v_1, \ldots, v_k)$, there exists a unique probability measure $\mathbb{Q}_\pi$ on $\Delta$, such that for all measurable $A \subseteq \overline{V}^{\mathbb{N}_0}$, we have*

$$(2.13) \quad P_\mathbf{0}((X_{k+t})_{t\in\mathbb{N}_0} \in A | (X_s)_{s=0,\ldots,k} = \pi) = \int_\Delta Q_{v_k,x}(A)\,\mathbb{Q}_\pi(dx).$$

*The measure $\mathbb{Q}_\pi$ is absolutely continuous with respect to $\mathbb{Q}$; the Radon–Nikodym derivative is given by*

$$(2.14) \quad \frac{d\mathbb{Q}_\pi}{d\mathbb{Q}}(x) = \frac{Q_{\mathbf{0},x}((X_s)_{s=0,\ldots,k} = \pi)}{P_\mathbf{0}((X_s)_{s=0,\ldots,k} = \pi)}, \qquad x \in \Delta.$$

The next theorem is concerned with the infinite-volume limit of the random environment. It shows that the random environment for finite ladders converges to the random environment for the infinite ladder as the size of the ladder grows.

THEOREM 2.8 (Infinite-volume limit). *As $n \to \infty$, the finite-dimensional marginals of $\mathbb{Q}^{(n)}$ converge weakly to the corresponding marginals of $\mathbb{Q}$.*

The next theorem tells us that ratios of the random weights of edges do not fluctuate too much: For neighboring edges $e$, $e'$, the random variables $x_e/x_{e'}$ are tight. The theorem makes even a more quantitative statement:

THEOREM 2.9 (Tail behavior of edge weights). *Let $e, e'$ be edges on level $i, j$, respectively. With respect to $\mathbb{Q}$, the random variable*

$$(2.15) \quad \ln\frac{x_e}{x_{e'}}$$

*has exponential tails. More precisely, there exist positive constants $c_7(a), c_8(a)$ such that for all $i, j \in \mathbb{N}$, edges $e, e'$ as above, and $M > 0$, one has*

$$(2.16) \quad \mathbb{Q}\left[\left|\ln\frac{x_e}{x_{e'}}\right| \geq M\right] \leq c_7 \max\{|i-j|, 1\} \exp\left\{-\frac{c_8 M}{\max\{|i-j|, 1\}}\right\}.$$

2.3. *Convergence to equilibrium.* In this subsection, we extend and refine the convergence result of Theorem 2.3. First, we describe the equilibrium measures $\mu_{\text{even}}$ and $\mu_{\text{odd}}$ from that theorem in terms of the random environment.

We define $x_{\text{even}}, x_{\text{odd}} : \overline{V} \to [0,1]$, by

$$(2.17) \quad x_{\text{even}}(v) := 1_{\overline{V}_{\text{even}}}(v) x_v \quad \text{and} \quad x_{\text{odd}}(v) := 1_{\overline{V}_{\text{odd}}}(v) x_v.$$

Since every $e \in \overline{E}$ contains one vertex in $\overline{V}_{\text{even}}$ and one vertex in $\overline{V}_{\text{odd}}$, we have

$$(2.18) \quad \sum_{v\in\overline{V}} x_{\text{even}}(v) = \sum_{v\in\overline{V}} x_{\text{odd}}(v) = \sum_{e\in\overline{E}} x_e = 1.$$

Hence $x_{\text{even}}$ and $x_{\text{odd}}$ are probability functions. By averaging these probability functions over the environment, we get the equilibrium distributions $\mu_{\text{even}}$ and $\mu_{\text{odd}}$. This is the claim of the following theorem:

THEOREM 2.10 (Annealed equilibrium measures are mixtures). *The measures $\mu_{\text{even}}$ and $\mu_{\text{odd}}$ from Theorem 2.3 have the following representation:*

$$\mu_{\text{even}}(v) = \int_\Delta x_{\text{even}}(v)\,\mathbb{Q}(dx),$$
(2.19)
$$\mu_{\text{odd}}(v) = \int_\Delta x_{\text{odd}}(v)\,\mathbb{Q}(dx) \qquad (v \in \overline{V}).$$

For any probability function $\mu: \overline{V} \to [0,1]$ and $x \in \Delta$, we define the law of the random walk in the environment $x$ with starting distribution $\mu$:

(2.20) $$Q_{\mu,x} := \sum_{v \in \overline{V}} \mu(v) Q_{v,x}.$$

In analogy to the convergence in distribution of $X_{2t}$ and $X_{2t+1}$ as stated in Theorem 2.3, there is a convergence result for the whole process:

THEOREM 2.11 (Convergence to equilibrium of the process). *As $s \to \infty$, the distributions of $(X_{2s+t})_{t \in \mathbb{N}_0}$ and $(X_{2s+1+t})_{t \in \mathbb{N}_0}$ converge in the following sense: For all measurable sets $A \subseteq \overline{V}^{\mathbb{N}_0}$, one has*

(2.21) $$P_{\mathbf{0}}((X_{2s+t})_{t \in \mathbb{N}_0} \in A) \stackrel{s \to \infty}{\longrightarrow} \int_\Delta Q_{x_{\text{even}},x}(A)\,\mathbb{Q}(dx),$$

(2.22) $$P_{\mathbf{0}}((X_{2s+1+t})_{t \in \mathbb{N}_0} \in A) \stackrel{s \to \infty}{\longrightarrow} \int_\Delta Q_{x_{\text{odd}},x}(A)\,\mathbb{Q}(dx).$$

We can interpret $x_{\text{even}}(v)$ and $x_{\text{odd}}(v)$ as the asymptotic probabilities to visit a vertex $v$ at even and at odd times, respectively, as time goes to infinity, conditioned on the environment. Even more, conditioned on the environment, the whole shifted process converges to a Markov chain, started in equilibrium:

THEOREM 2.12 (Convergence to equilibrium, conditioned on the environment). *For all $x \in \Delta$, the laws of $X_{2t}$ and of $X_{2t+1}$ with respect to $Q_{\mathbf{0},x}$ converge:*

(2.23) $$\lim_{t \to \infty} Q_{\mathbf{0},x}(X_{2t} = v) = x_{\text{even}}(v),$$

(2.24) $$\lim_{t \to \infty} Q_{\mathbf{0},x}(X_{2t+1} = v) = x_{\text{odd}}(v) \qquad \text{for all } v \in \overline{V}.$$

*More generally, for all $A \subseteq \overline{V}^{\mathbb{N}_0}$ measurable,*

$$Q_{\mathbf{0},x}((X_{2s+t})_{t\in\mathbb{N}_0} \in A) \stackrel{s\to\infty}{\longrightarrow} Q_{x_{\mathrm{even}},x}(A) \tag{2.25}$$

*and*

$$Q_{\mathbf{0},x}((X_{2s+1+t})_{t\in\mathbb{N}_0} \in A) \stackrel{s\to\infty}{\longrightarrow} Q_{x_{\mathrm{odd}},x}(A). \tag{2.26}$$

2.4. *A representation of the random environment by an infinite-volume Gibbs measure.* In this subsection, we show that the random environment for the reinforced random walk on the *infinite* graph $\overline{G}$ can be written in terms of an *infinite-volume* Gibbs measure. For finite pieces of $\overline{G}$, the description in terms of *finite-volume* Gibbs measures is one of the central ideas in [7] and [10]. Here, we deal with the thermodynamic limit of these Gibbs measures. But first, we review the state spaces for the local spin variables; for more details in the finite-dimensional setup see [7] and [10]. Roughly speaking, every level of the ladder corresponds to a "compound spin variable." Although the precise form of the state spaces is irrelevant for almost all arguments that will follow, we explain very briefly their intuitive meaning: The rungs of the ladder are the subgraphs $\{i\} \times G$, that is, the vertical edges if we view the ladder as going from the left to the right. On the other hand, the (horizontal) edges connecting these subgraphs constitute the "slices" of the ladder.

Very roughly speaking, the spin variables consist of logarithms of ratios of neighboring edge weights. Besides, the spin variables have also discrete components. These discrete components provide local descriptions of spanning trees on the one hand and on the other hand some signs which are not encoded in the logarithms above. The precise connection between the spin variables and the edge weights is given in Definitions 2.19 and 3.2, below.

In the whole article, we use the more general notation, the state spaces, and the variable transforms from [10] rather than [7]. However, in the special case $\overline{G} = \mathbb{N}_0 \times \{1,2\}$, the variant described in [7] could also be used. This variant uses a slightly different variable transform, and it applies to all initial weights $a > 3/4$ rather than only large $a$. Therefore, we include below also citations of [7]. Readers interested only in large $a$ may ignore them.

DEFINITION 2.13 (*State spaces*). As in Definition 2.15 of [10] (see also Definition 2.8 of [7]), we fix two different vertices $v^{\mathrm{tree}}, v^* \in V$, and we set

$$\Omega_{\mathrm{left}} := \mathbb{R}^E, \tag{2.27}$$

$$\Omega_{\mathrm{slice}} := \mathbb{R}^V \times \{\pm 1\}^{V \setminus \{v^{\mathrm{tree}},v^*\}} \times \mathbb{R}^{V \setminus \{v^{\mathrm{tree}},v^*\}} \times \mathrm{Treevar},$$

$$\Omega_{\mathrm{rung}} := \mathbb{R}^E \times \mathbb{R}; \tag{2.28}$$

here Treevar is a finite set (see Definition 2.7 of [10]). Furthermore, we define the index set

$$\mathcal{I} := \{\text{left}\} \cup \{(\text{rung}, i), (\text{slice}, i) \colon i \in \mathbb{N}\}. \tag{2.29}$$

For any $\Lambda \subseteq \mathcal{I}$, we define

$$\Omega_\Lambda := \prod_{\iota \in \Lambda} \Omega_\iota, \tag{2.30}$$

where $\Omega_{\text{rung},i} = \Omega_{\text{rung}}$ and $\Omega_{\text{slice},i} = \Omega_{\text{slice}}$ for $i \in \mathbb{N}$. We abbreviate $\Omega = \Omega_\mathcal{I}$. The $\sigma$-algebra on $\Omega_\Lambda$, induced by the canonical projections, is denoted by $\mathcal{B}(\Omega_\Lambda)$.

For most parts of our arguments, even the precise form of the components of the compound spin variables is irrelevant. Thus, we introduce abbreviations for the compound spins:

DEFINITION 2.14 (*Local state variables*). For $\Lambda \subseteq \mathcal{I}$, we denote the canonical element of $\Omega_\Lambda$ by $\omega_\Lambda := (\omega_\iota)_{\iota \in \Lambda}$. We use the following notation from Definition 2.15 of [10] for the components of $\omega_\iota$ (see also Definition 2.8 of [7]):

$$\omega_{\text{slice},i} := ((X_i(v))_{v \in V}, (\sigma_i(v), W_i(v))_{v \in V \setminus \{v^{\text{tree}}, v^*\}}, T_i), \tag{2.31}$$

$$\omega_{\text{rung},i} := ((Z_i(e))_{e \in E}, \Gamma_i), \tag{2.32}$$

where $i \in \mathbb{N}$. In addition,

$$\omega_{\text{left}} := ((Z_0(e))_{e \in E}). \tag{2.33}$$

The components $X_i(v)$, $W_i(v)$, $Z_i(e)$, and $\Gamma_i$ take values in $\mathbb{R}$, whereas $\sigma_i(v)$ and $T_i$ take only finitely many values. We denote by $d\omega_\Lambda$ the Lebesgue measure on the continuous components times the counting measure on the discrete components.

For $n \in \mathbb{N}$, we consider the index set

$$[0, n] := \{\text{left}\} \cup \{(\text{rung}, i), (\text{slice}, i) \colon i = 1, \ldots, n\} \subset \mathcal{I} \tag{2.34}$$

with boundary $\partial[0, n] := \{(\text{slice}, n+1)\}$.

Let $H_{\text{middle}} := H_{\text{middle},a,1/4}$ and $H_{\text{left}} := H_{\text{left},a}$ be the local Hamiltonians as defined in Definitions 2.18 and 2.20, respectively, of [10] (see also Definitions 2.10 and 2.11 of [7]). Their explicit form is irrelevant for this paper. Since the initial weight $a$ is kept constant, we suppress it in the notation. Using these local Hamiltonians, the finite-volume Gibbs measures are defined in the standard way as follows.

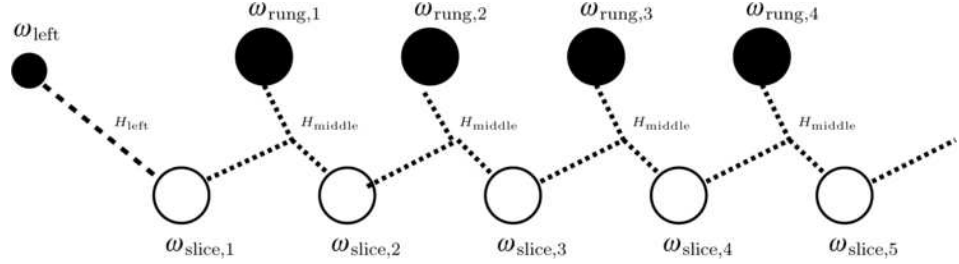

Fig. 1. *Interactions in the Gibbs measure.*

DEFINITION 2.15 (*Finite-volume Gibbs measures with boundary conditions*). Let $n \in \mathbb{N}$ and $\Lambda = [0, n]$. We define the finite-volume Hamiltonian with boundary conditions $\omega_{\partial \Lambda}$ by

$$
\begin{aligned}
(2.35) \quad H_{[0,n]}(\omega_{[0,n]} | \omega_{\partial[0,n]}) &:= H_{\text{left}}(\omega_{\text{left}}, \omega_{\text{slice},1}) \\
&\quad + \sum_{i=1}^{n} H_{\text{middle}}(\omega_{\text{slice},i}, \omega_{\text{rung},i}, \omega_{\text{slice},i+1}).
\end{aligned}
$$

Furthermore, for any bounded and measurable function $F : \Omega_\Lambda \to \mathbb{R}$, we introduce $g_\Lambda^F : \Omega_{\partial \Lambda} \to \mathbb{R}$ by

$$
(2.36) \qquad g_\Lambda^F(\omega_{\partial \Lambda}) := \int_{\Omega_\Lambda} F(\omega_\Lambda) e^{-H_\Lambda(\omega_\Lambda | \omega_{\partial \Lambda})} \, d\omega_\Lambda.
$$

The finite-volume Gibbs measure $\mathbb{P}_\Lambda(\cdot | \omega_{\partial \Lambda})$ is given by the Markov kernel $\mathbb{K}_\Lambda : \mathcal{B}(\Omega_\Lambda) \times \Omega_{\partial \Lambda} \to [0, 1]$, where

$$
(2.37) \qquad \mathbb{K}_\Lambda(A, \omega_{\partial \Lambda}) := \frac{g_\Lambda^{1_A}(\omega_{\partial \Lambda})}{g_\Lambda^1(\omega_{\partial \Lambda})}.
$$

Figure 1 illustrates the interactions in the Gibbs measures as described by the local Hamiltonians. The local interactions $H_{\text{middle}}$ are drawn nonsymmetrically; this represents symbolically the absence of reflection symmetry caused by the symmetry-breaking terms $\Gamma_i/4$ arising in the summand $H_{\text{middle}}(\omega_{\text{slice},i}, \omega_{\text{rung},i}, \omega_{\text{slice},i+1})$ (see (2.41) of [10] and (2.42) of [7]; recall that in the notation of these articles, $H_{\text{middle}} = H_{\text{middle},a,1/4}$).

A central role in our analysis is played by the transfer operator and its leading eigenfunctions from the left and from the right. The transfer operator is an integral operator with an integral kernel $k$ defined as follows: As in Definition 4.1 and in particular (4.4) of [10] (see also Definition 4.2 and (4.5) of [7]), we define

$$
(2.38) \qquad k(\omega_{\text{slice}}, \omega'_{\text{slice}}) := \int_{\Omega_{\text{rung}}} e^{-H_{\text{middle}}(\omega_{\text{slice}}, \omega_{\text{rung}}, \omega'_{\text{slice}})} \, d\omega_{\text{rung}};
$$

in the notation of [7] and [10], $k = k_\eta^\Upsilon$ for $\eta = 1/4$ and $\Upsilon \equiv 1$. The transfer operators, that is, the (left and right) integral operators with integral kernel $k$, have a positive leading eigenvalue $\lambda$, as was shown in Lemma 4.2 of [10] (see also Lemma 4.4 of [7]). Let $v$, $v^*$ denote their leading eigenfunctions; they are positive and unique up to normalization:

$$(2.39) \qquad \int_{\Omega_{\text{slice}}} k(\omega_{\text{slice}}, \omega'_{\text{slice}}) v^*(\omega'_{\text{slice}}) \, d\omega'_{\text{slice}} = \lambda v^*(\omega_{\text{slice}}),$$

$$(2.40) \qquad \int_{\Omega_{\text{slice}}} v(\omega_{\text{slice}}) k(\omega_{\text{slice}}, \omega'_{\text{slice}}) \, d\omega_{\text{slice}} = \lambda v(\omega'_{\text{slice}}).$$

The infinite-volume Gibbs measures are defined by averaging finite volume Gibbs measures with boundary conditions weighted according to the right eigenfunction of the transfer operator, just as always in one-dimensional Gibbs systems with short-range interactions.

DEFINITION 2.16 (*Infinite-volume Gibbs measure*). We define the infinite-volume Gibbs measure $\mathbb{P}$ as the unique probability measure on $\mathcal{B}(\Omega)$ satisfying

$$(2.41) \qquad \mathbb{P}(\omega_\Lambda \in A) = \frac{\int_{\Omega_{\partial\Lambda}} g_\Lambda^{1_A}(\omega_{\partial\Lambda}) v^*(\omega_{\partial\Lambda}) \, d\omega_{\partial\Lambda}}{\int_{\Omega_{\partial\Lambda}} g_\Lambda^1(\omega_{\partial\Lambda}) v^*(\omega_{\partial\Lambda}) \, d\omega_{\partial\Lambda}}$$

for all $\Lambda = [0, n]$ and $A \in \mathcal{B}(\Omega_\Lambda)$.

LEMMA 2.17 (Infinite-volume Gibbs measure). *The infinite-volume Gibbs measure $\mathbb{P}$ is well defined.*

Indeed, the infinite-volume Gibbs measure satisfies the DLR-conditions:

THEOREM 2.18 (Dobrushin–Landford–Ruelle conditions). *The kernel $\mathbb{K}_\Lambda$ is a regular conditional distribution of $\omega_\Lambda$ with respect to $\mathbb{P}$ conditioned on $\omega_{\partial\Lambda}$. In particular, for any set $A \in \mathcal{B}(\Omega_\Lambda)$, we have*

$$(2.42) \qquad \mathbb{P}(\omega_\Lambda \in A | \omega_{\partial\Lambda}) = \mathbb{K}_\Lambda(A, \omega_{\partial\Lambda}).$$

*Slightly more generally, for any bounded measurable function $F : \Omega \to \mathbb{R}$*

$$(2.43) \qquad \int_\Omega F \, d\mathbb{P} = \int_\Omega \int_{\Omega_\Lambda} F(\chi_\Lambda, \omega_{\Lambda^c}) \mathbb{K}_\Lambda(d\chi_\Lambda, \omega_{\partial\Lambda}) \, \mathbb{P}(d\omega).$$

Next, we define a transformation from $\omega \in \Omega$ to the edge weights $(\tilde{x}_e)_{e \in \overline{E}} \in \mathbb{R}_+^{\overline{E}}$. Up to a normalization, the $\tilde{x}_e$ turn out to be the random environment for the reinforced random walk; see Theorem 2.21, below. The transformation is the same as in Definition 2.17 of [10] with the abbreviations from Definition 2.16 of [10] plugged in (see also Definition 2.9 of [7]). The form of this transformation is motivated in [7] and [10]. Its precise form is used in [7] and [10] to derive bounds which are also relevant in the present paper. However, the precise form of the transformation is not used below.

For the following definition, recall the notation $e_i$ for edges, introduced after formula (1.4).

DEFINITION 2.19 (*Transformation from state variables to random environment*). For $\omega \in \Omega$, we define

$$\tilde{x}(\omega) = (\tilde{x}_e(\omega))_{e \in \overline{E}} \in \mathbb{R}_+^{\overline{E}} = (0, \infty)^{\overline{E}} \tag{2.44}$$

as follows: We set for $v \in V$ and $i \in \mathbb{N}$

$$\tilde{x}_{\{v_{i-1}, v_i\}} := e^{X_1(v^*) - Z_0(e^*)} \exp\left[X_i(v) - X_i(v^*) - \sum_{j=1}^{i-1} \Gamma_j\right] \tag{2.45}$$

and for $e \in E$

$$\tilde{x}_{e_0} := \exp[Z_0(e) - Z_0(e^*)], \tag{2.46}$$

$$\tilde{x}_{e_i} := e^{X_1(v^*) - Z_0(e^*)} \tag{2.47}$$
$$\times \exp\left[Z_i(e) - \tfrac{1}{2}\{X_i(v^*) + X_{i+1}(v^*) + \Gamma_i\} - \sum_{j=1}^{i-1} \Gamma_j\right].$$

We define $\tilde{\mathbb{Q}}$ to be the law of $\tilde{x}$ with respect to $\mathbb{P}$.

Clearly, $\tilde{x}_{e_0^*} = 1$. We want to change the normalization of the edge weights $(\tilde{x}_e)_e$ in such a way that they sum up to one. This is done in the following definition. Note that multiplying all edge weights $(x_e)_{e \in \overline{E}}$ by the same positive constant does not change the measure $Q_{v,x}$.

DEFINITION 2.20 (*Changing the normalization*). We define $x(\omega) = (x_e(\omega))_{e \in \overline{E}} \in \Delta$ by

$$x_e := \frac{\tilde{x}_e}{\sum_{e' \in \overline{E}} \tilde{x}_{e'}}, \tag{2.48}$$

whenever these random variables are well defined.

By a slight abuse of notation, we use the same symbol $x_e$ in two slightly different ways: On the one hand, $x = (x_e)_{e \in \overline{E}}$ denotes weights of the random environment, for example, in Theorem 2.4. On the other hand, for example, in (2.48), $x_e(\omega)$ denotes the value of a random variable on $\Omega$. The following theorem justifies this little abuse of notation:

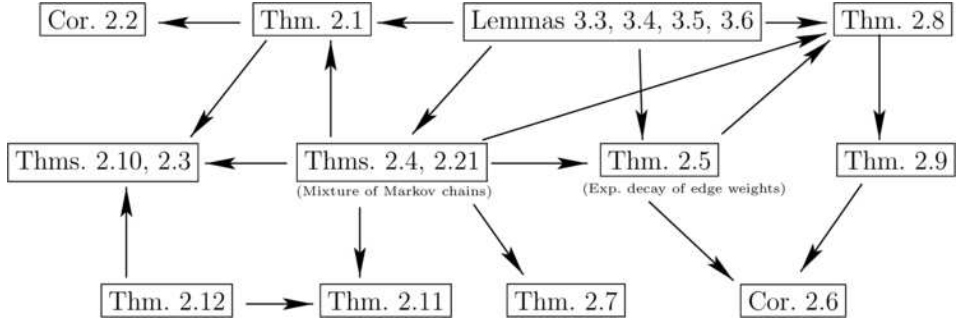

Fig. 2.

Theorem 2.21 (Representation by an infinite-volume Gibbs measure). *Let $\omega$ be a random variable with distribution $\mathbb{P}$. Then, the distribution $\tilde{\mathbb{Q}}$ of $\tilde{x}(\omega)$ equals the distribution of the random environment, normalized such that the reference edge $e_0^*$ gets weight $\tilde{x}_{e_0^*} = 1$.*

*In addition, $x(\omega)$ is almost surely well defined. Its distribution equals $\mathbb{Q}$, the distribution of the random environment, normalized such that $\sum_{e \in \overline{E}} x_e = 1$.*

## 3. Proofs.

*Organization of this paper.* The order in which the results are proven differs from the order in which we presented them. The complicated dependence structure between the theorems is best represented graphically: Figure 2 displays the mutual dependence of the lemmas, theorems and corollaries in this paper. Lemma 2.17 ($\mathbb{P}$ is well defined) and Theorem 2.18 (Dobrushin–Landford–Ruelle conditions) are not represented in Figure 2.

3.1. *Analysis of the infinite-volume Gibbs measure.*

Proof of Lemma 2.17. First, observe that the integrands arising in (2.41) decay exponentially at infinity. For $H_{\text{middle}}$ and $H_{\text{left}}$, this is shown in Propositions 3.1 and 3.9 of [10] (see also Propositions 3.2 and 3.5 of [7]). Because of the eigenvalue equation (2.39), the decay of the local Hamiltonians implies the exponential decay of $\upsilon^*$. Thus, the numerator and the denominator of the fraction in (2.41) are both finite. The denominator is strictly positive.

The set of events of the form "$\omega_\Lambda \in A$" is closed under intersections and generates $\mathcal{B}(\Omega)$. Hence, it remains to show that the definition (2.41) does not depend on the choice of $\Lambda$, that is,

$$(3.1) \qquad \mathbb{P}(\omega_{\Lambda'} \in A') = \mathbb{P}(\omega_\Lambda \in A)$$

holds whenever $\Lambda \subset \Lambda'$, $A \in \mathcal{B}(\Omega_\Lambda)$, and $A'$ is the inverse image of $A$ with respect to the projection $\Omega_{\Lambda'} \to \Omega_\Lambda$. It suffices to take $\Lambda'$ one step bigger than $\Lambda$.

Let $\Lambda = [0, n-1]$, $\Lambda' = [0, n]$. By the definition (2.41) of $\mathbb{P}$,

$$(3.2) \quad \mathbb{P}(\omega_{\Lambda'} \in A') := \frac{\int_{\Omega_{\partial\Lambda'}} g_{\Lambda'}^{1_{A'}}(\omega_{\partial\Lambda'}) v^*(\omega_{\partial\Lambda'}) \, d\omega_{\partial\Lambda'}}{\int_{\Omega_{\partial\Lambda'}} g_{\Lambda'}^{1}(\omega_{\partial\Lambda'}) v^*(\omega_{\partial\Lambda'}) \, d\omega_{\partial\Lambda'}}.$$

Observe that $\Lambda' = \Lambda \cup \partial\Lambda \cup \{(\text{rung}, n)\}$ and

$$(3.3) \quad H_{\Lambda'}(\omega_{\Lambda'} | \omega_{\partial\Lambda'}) = H_\Lambda(\omega_\Lambda | \omega_{\partial\Lambda}) + H_{\text{middle}}(\omega_{\partial\Lambda}, \omega_{\text{rung},n}, \omega_{\partial\Lambda'}).$$

Hence, by Fubini's theorem,

$$(3.4) \quad \begin{aligned} g_{\Lambda'}^{1_{A'}}(\omega_{\partial\Lambda'}) &= \int_{\Omega_{\Lambda'}} d\omega_{\Lambda'} 1_{A'}(\omega_{\Lambda'}) e^{-H_{\Lambda'}(\omega_{\Lambda'} | \omega_{\partial\Lambda'})} \\ &= \int_{\Omega_{\partial\Lambda}} d\omega_{\partial\Lambda} \int_{\Omega_{\text{rung}}} d\omega_{\text{rung}} e^{-H_{\text{middle}}(\omega_{\partial\Lambda}, \omega_{\text{rung}}, \omega_{\partial\Lambda'})} \\ &\quad \times \int_{\Omega_\Lambda} d\omega_\Lambda 1_A(\omega_\Lambda) e^{-H_\Lambda(\omega_\Lambda | \omega_{\partial\Lambda})} \\ &= \int_{\Omega_{\partial\Lambda}} d\omega_{\partial\Lambda} \, k(\omega_{\partial\Lambda}, \omega_{\partial\Lambda'}) g^{1_A}(\omega_{\partial\Lambda}). \end{aligned}$$

Rewriting (2.39) as

$$(3.5) \quad \int_{\Omega_{\partial\Lambda'}} k(\omega_{\partial\Lambda}, \omega_{\partial\Lambda'}) v^*(\omega_{\partial\Lambda'}) \, d\omega_{\partial\Lambda'} = \lambda v^*(\omega_{\partial\Lambda})$$

and applying Fubini's theorem, we conclude that

$$(3.6) \quad \begin{aligned} \int_{\Omega_{\partial\Lambda'}} &g_{\Lambda'}^{1_{A'}}(\omega_{\partial\Lambda'}) v^*(\omega_{\partial\Lambda'}) \, d\omega_{\partial\Lambda'} \\ &= \lambda \int_{\Omega_{\partial\Lambda}} g_\Lambda^{1_A}(\omega_{\partial\Lambda}) v^*(\omega_{\partial\Lambda}) \, d\omega_{\partial\Lambda}. \end{aligned}$$

The last identity with $A$ replaced by $\Omega_\Lambda$ and $A'$ replaced by $\Omega_{\Lambda'}$ gives an identity for the denominator in (3.2). Combining these two identities yields the claim (3.1). $\square$

PROOF OF THEOREM 2.18. It suffices to prove (2.43) for nonnegative functions $F$ that depend on only finitely many components, say $F:\Omega_{\Lambda'} \to \mathbb{R}$, where $\Lambda' \supseteq \Lambda \cup \partial\Lambda$. We decompose $H_{\Lambda'}$ as

$$(3.7) \quad H_{\Lambda'}(\omega_{\Lambda'} | \omega_{\partial\Lambda'}) = H_\Lambda(\omega_\Lambda | \omega_{\partial\Lambda}) + H_{\Lambda',\Lambda}(\omega_{\Lambda' \setminus \Lambda} | \omega_{\partial\Lambda'})$$

with a sum of local Hamiltonians $H_{\Lambda',\Lambda}(\omega_{\Lambda'\setminus\Lambda}|\omega_{\partial\Lambda'}) := H_{\Lambda'}(\omega_{\Lambda'}|\omega_{\partial\Lambda'}) - H_{\Lambda}(\omega_{\Lambda}|\omega_{\partial\Lambda})$ that depends only on $\omega_{\Lambda'\setminus\Lambda}$ and $\omega_{\partial\Lambda'}$. We calculate for $\omega_{\partial\Lambda'} \in \Omega_{\partial\Lambda'}$:

$$
\begin{aligned}
\int_{\Omega_{\Lambda'}} &d\omega_{\Lambda'} e^{-H_{\Lambda'}(\omega_{\Lambda'}|\omega_{\partial\Lambda'})} F(\omega_{\Lambda'}) \\
&= \int_{\Omega_\Lambda} d\chi_\Lambda \int_{\Omega_{\Lambda'\setminus\Lambda}} d\omega_{\Lambda'\setminus\Lambda} e^{-H_\Lambda(\chi_\Lambda|\omega_{\partial\Lambda})} e^{-H_{\Lambda',\Lambda}(\omega_{\Lambda'\setminus\Lambda}|\omega_{\partial\Lambda'})} \\
&\quad \times \frac{F(\chi_\Lambda, \omega_{\Lambda'\setminus\Lambda})}{g^1_\Lambda(\omega_{\partial\Lambda})} \int_{\Omega_\Lambda} d\omega_\Lambda e^{-H_\Lambda(\omega_\Lambda|\omega_{\partial\Lambda})} \\
(3.8)\quad &= \int_{\Omega_\Lambda} d\omega_\Lambda \int_{\Omega_{\Lambda'\setminus\Lambda}} d\omega_{\Lambda'\setminus\Lambda} e^{-H_\Lambda(\omega_\Lambda|\omega_{\partial\Lambda})} e^{-H_{\Lambda',\Lambda}(\omega_{\Lambda'\setminus\Lambda}|\omega_{\partial\Lambda'})} \\
&\quad \times \int_{\Omega_\Lambda} d\chi_\Lambda e^{-H_\Lambda(\chi_\Lambda|\omega_{\partial\Lambda})} \frac{F(\chi_\Lambda, \omega_{\Lambda'\setminus\Lambda})}{g^1_\Lambda(\omega_{\partial\Lambda})} \\
&= \int_{\Omega_{\Lambda'}} d\omega_{\Lambda'} e^{-H_{\Lambda'}(\omega_{\Lambda'}|\omega_{\partial\Lambda'})} \int_{\Omega_\Lambda} \mathbb{K}_\Lambda(d\chi_\Lambda, \omega_{\partial\Lambda}) F(\chi_\Lambda, \omega_{\Lambda'\setminus\Lambda});
\end{aligned}
$$

recall definition (2.37) of the Markov kernel $\mathbb{K}_\Lambda$. Claim (2.43) follows using the definition (2.41) of $\mathbb{P}$: one multiplies the finite-volume DLR-equation (3.8) by the eigenfunction $v^*(\omega_{\partial\Lambda'})$, integrates with respect to $d\omega_{\partial\Lambda'}$, and divides by the normalizing constant $\int_{\Omega_{\partial\Lambda'}} g^1_{\Lambda'}(\omega_{\partial\Lambda'}) v^*(\omega_{\partial\Lambda'}) d\omega_{\partial\Lambda'}$. $\square$

Recall that $\overline{G}^{(n)} = (\overline{V}^{(n)}, \overline{E}^{(n)})$ denotes the restriction of the graph $\overline{G}$ to the finite piece $\overline{V}^{(n)} = \{0, \ldots, n\} \times V$.

DEFINITION 3.1 (*Finite-volume Gibbs measures*). For $n \in \mathbb{N}$, we define in analogy to Definition 2.13:

(3.9) $\Omega_{\text{right}} := \mathbb{R}^E$,

$\mathcal{I}^{(n)} := \{\text{left}\} \cup \{(\text{rung}, i), (\text{slice}, i) : 1 \leq i \leq n-1\} \cup \{(\text{slice}, n), \text{right}\}$

(3.10) $\quad = [0, n-1] \cup \{(\text{slice}, n), \text{right}\}$,

(3.11) $\Omega^{(n)} := \prod_{\iota \in \mathcal{I}^{(n)}} \Omega_\iota$.

The canonical element of $\Omega^{(n)}$ is denoted by $\omega^{(n)}$.

Let $H_{\text{right}} := H_{\text{right},a}$ be as in Definition 2.21 of [10] (see also Definition 2.11 of [7]). We define the finite-volume Hamiltonian over $\mathcal{I}^{(n)}$ as follows:

$$H^{(n)}(\omega^{(n)}) := H_{\text{left}}(\omega_{\text{left}}, \omega_{\text{slice},1})$$

$$+ \sum_{i=1}^{n-1} H_{\text{middle}}(\omega_{\text{slice},i}, \omega_{\text{rung},i}, \omega_{\text{slice},i+1})$$

(3.12)
$$+ H_{\text{right}}(\omega_{\text{slice},n}, \omega_{\text{right}})$$
$$= H_{[0,n-1]}(\omega_{[0,n-1]}|\omega_{\partial[0,n-1]}) + H_{\text{right}}(\omega_{\partial[0,n-1]}, \omega_{\text{right}}).$$

Finally, the finite-volume Gibbs measure $\mathbb{P}^{(n)}$ is defined to be the probability measure given by

(3.13) $$\mathbb{P}^{(n)}(\omega^{(n)} \in A) := \frac{\int_{\Omega^{(n)}} 1_A(\omega^{(n)}) e^{-H^{(n)}(\omega^{(n)})} \, d\omega^{(n)}}{\int_{\Omega^{(n)}} e^{-H^{(n)}(\omega^{(n)})} \, d\omega^{(n)}}.$$

For $A \in \mathcal{B}(\Omega_{[0,n-1]})$, this comes down to

(3.14)
$$\mathbb{P}^{(n)}(\omega_{[0,n-1]} \in A) = \frac{\int_{\Omega_{\text{slice}}} g^{1_A}_{[0,n-1]}(\omega_{\text{slice}}) g_{\text{right}}(\omega_{\text{slice}}) \, d\omega_{\text{slice}}}{\int_{\Omega_{\text{slice}}} g^{1}_{[0,n-1]}(\omega_{\text{slice}}) g_{\text{right}}(\omega_{\text{slice}}) \, d\omega_{\text{slice}}}$$
$$= \frac{\langle g^{1_A}_{[0,n-1]} g_{\text{right}} \rangle}{\langle g^{1}_{[0,n-1]} g_{\text{right}} \rangle},$$

where

(3.15) $$g_{\text{right}}(\omega_{\text{slice}}) = \int_{\Omega_{\text{right}}} e^{-H_{\text{right}}(\omega_{\text{slice}}, \omega_{\text{right}})} \, d\omega_{\text{right}}$$

and $\langle fg \rangle$ is a short notation for $\int_{\Omega_{\text{slice}}} f(\omega_{\text{slice}}) g(\omega_{\text{slice}}) \, d\omega_{\text{slice}}$.

DEFINITION 3.2 (*Transformation of variables—finite volume version*). For $\omega^{(n)} \in \Omega^{(n)}$, we define

(3.16) $$\tilde{x}^{(n)}(\omega^{(n)}) = (\tilde{x}^{(n)}_e(\omega^{(n)}))_{e \in \overline{E}^{(n)}} \in \mathbb{R}_+^{\overline{E}^{(n)}}$$

as follows: Whenever $e$ is not an edge in the right border of the finite ladder $\overline{G}^{(n)}$, that is, whenever $e \neq \{u_n, v_n\}$ for all $\{u, v\} \in E$, we define $\tilde{x}^{(n)}_e$ just as in (2.45), (2.46) and (2.47) in Definition 2.19. However, if $e$ is an edge in the right border of the finite ladder $\overline{G}^{(n)}$, that is, if $e = e_n = \{u_n, v_n\}$, we set

(3.17) $$\tilde{x}^{(n)}_{e_n} := e^{X_1(v^*) - Z_0(e^*)} \exp\left[Z_n(e) - X_n(v^*) - \sum_{j=1}^{n-1} \Gamma_j\right].$$

Let $\tilde{\mathbb{Q}}^{(n)}$ denote the law of $\tilde{x}^{(n)}(\omega^{(n)}) = (\tilde{x}^{(n)}_e(\omega^{(n)}))_{e \in \overline{E}^{(n)}}$, provided that $\omega^{(n)}$ has the distribution $\mathbb{P}^{(n)}$. Finally, we set

(3.18) $$x^{(n)}_e := \frac{\tilde{x}^{(n)}_e}{\sum_{e' \in \overline{E}^{(n)}} \tilde{x}^{(n)}_{e'}}.$$

Table 1

| Random variable | Normalization | Law |
|---|---|---|
| Infinite-volume | | |
| State variable $\omega$ | | $\mathbb{P}$ |
| Random environment $\tilde{x}$ | $\tilde{x}_{e_0^*} = 1$ | $\tilde{\mathbb{Q}}$ |
| Random environment $x$ | $\sum_{e \in \overline{E}} x_e = 1$ | $\mathbb{Q}$ |
| Finite-volume | | |
| State variable $\omega^{(n)}$ | | $\mathbb{P}^{(n)}$ |
| Random environment $\tilde{x}^{(n)}$ | $\tilde{x}_{e_0^*}^{(n)} = 1$ | $\tilde{\mathbb{Q}}^{(n)}$ |
| Random environment $x^{(n)}$ | $\sum_{e \in \overline{E}^{(n)}} x_e^{(n)} = 1$ | $\mathbb{Q}^{(n)}$ |

Thus, we use two different normalizations: The "tilde" version, where the reference edge $e_0^*$ gets weight 1, and the "no-tilde" version, where all weights sum up to 1. Note that, in addition to (3.18), we have the following conversion between the "tilde" and "no-tilde" normalization:

$$(3.19) \qquad \tilde{x}_e^{(n)} = \frac{x_e^{(n)}}{x_{e_0^*}^{(n)}}, \qquad \tilde{x}_e = \frac{x_e}{x_{e_0^*}} \quad \text{and} \quad x_e = \frac{\tilde{x}_e}{\sum_{e' \in \overline{E}} \tilde{x}_{e'}}.$$

It turns out that the $x_e^{(n)}(\omega^{(n)})$ have the distribution $\mathbb{Q}^{(n)}$, provided that $\omega^{(n)}$ has the law $\mathbb{P}^{(n)}$. This is a consequence of Lemma 3.3 below.

As a mnemonic aid, Table 1 summarizes the finite-volume and infinite-volume measures and the corresponding random variables.

Note that $\mathbb{P}$, $\tilde{\mathbb{Q}}$, $\mathbb{P}^{(n)}$, and $\tilde{\mathbb{Q}}^{(n)}$ are by definition the laws of the random variables $\omega$, $\tilde{x}(\omega)$, $\omega^{(n)}$ and $\tilde{x}^{(n)}(\omega^{(n)})$, respectively. However, the fact that $\mathbb{Q}$ is the law of $x(\omega)$ is stated in Theorem 2.21, and the fact that the $\mathbb{Q}^{(n)}$ is the law of $x^{(n)}(\omega^{(n)})$ is a consequence of the results in [7] and [10]. This is made precise in the following lemma:

Lemma 3.3 (Connection between Gibbs measure and random environment—finite-volume version). *The measure $\mathbb{Q}^{(n)}$, that is, the mixing measure describing the random environment for reinforced random walk on the finite ladder $\overline{G}^{(n)}$, equals the joint distribution of the random variables $(x_e^{(n)}(\omega^{(n)}))_{e \in \overline{E}^{(n)}}$, provided that $\omega^{(n)}$ has the law $\mathbb{P}^{(n)}$.*

*Moreover, let $F \subset \overline{E}$ be finite. Take $n$ large enough so that $F \subseteq \overline{E}^{(n-2)}$. Then $(\tilde{x}_e(\omega))_{e \in F}$ equals $(\tilde{x}_e^{(n)}(\omega^{(n)}))_{e \in F}$, provided that $\omega_{[0,n-1]} = \omega_{[0,n-1]}^{(n)}$.*

Proof. Combining Lemma 2.24 of [10] and Definition 2.17 of [10] (see also Lemma 2.13 of [7] and Definition 2.9 of [7]), the first claim follows.

To prove the second claim, note that each $\tilde{x}_e(\omega)$ depends only on *finitely many* components of $\omega$. More precisely, it depends only on $\omega_{[0,n-1]}$, if the edge $e$ is on a level strictly less than $n-1$. Consider a finite ladder $\overline{G}^{(n)}$ and any edge $e \in \overline{E}^{(n-2)}$. Then, the transformation $\omega \mapsto \tilde{x}_e(\omega)$ described for the *infinite* ladder in Definition 2.19 coincides with the map $\omega^{(n)} \mapsto \tilde{x}_e^{(n)}(\omega^{(n)})$ for the *finite* ladder $\overline{G}^{(n)}$, given in Definition 3.2; see also Definition 2.17 of [10] and Definition 2.9 of [7]. This proves the second claim. $\square$

LEMMA 3.4 (Thermodynamic limit). *As $n \to \infty$, the finite-dimensional marginals of $\mathbb{P}^{(n)}$ converge weakly to the corresponding marginals of $\mathbb{P}$. Even more, for any measurable bounded function $f$ depending only on finitely many coordinates $\omega_i$, $i \in I$, we have*

$$(3.20) \qquad \lim_{n \to \infty} \int f(\omega_i, i \in I) \mathbb{P}^{(n)}(d\omega) = \int f(\omega_i, i \in I) \mathbb{P}(d\omega).$$

*A similar statement holds for $\tilde{\mathbb{Q}}^{(n)} \to \tilde{\mathbb{Q}}$.*

PROOF. Let $K : L^2(\Omega_{\text{slice}}) \to L^2(\Omega_{\text{slice}})$,

$$(3.21) \qquad Kf(\omega_{\text{slice}}) := \int_{\Omega_{\text{slice}}} k(\omega_{\text{slice}}, \omega'_{\text{slice}}) f(\omega'_{\text{slice}}) \, d\omega'_{\text{slice}}$$

denote the integral operator with integral kernel $k$, defined in (2.38). Furthermore, normalizing the integral operator by its leading eigenvalue $\lambda$, we set $\hat{K} := \lambda^{-1} K$.

Let $l \in \mathbb{N}$, and take a bounded measurable function $f$ depending only on $\omega_{[0,l]}$. Take $n > l$. Then, using the boundary function $g_{\text{right}}$ from (3.15), we write

$$(3.22) \qquad \int f(\omega_\Lambda) \mathbb{P}^{(n)}(d\omega) = \frac{\langle g^f_{[0,l]} K^{n-l} g_{\text{right}} \rangle}{\langle g^1_{[0,l]} K^{n-l} g_{\text{right}} \rangle} = \frac{\langle g^f_{[0,l]} \hat{K}^{n-l} g_{\text{right}} \rangle}{\langle g^1_{[0,l]} \hat{K}^{n-l} g_{\text{right}} \rangle};$$

recall the definition (2.36) of $g^f_{[0,l]} \in L^2(\Omega_{\text{slice}})$. Now, for fixed $l \in \mathbb{N}$, Corollary 4.3 in [10] (see also Corollary 4.5 in [7]) states that

$$(3.23) \qquad \hat{K}^{n-l} g_{\text{right}} \stackrel{n \to \infty}{\longrightarrow} \upsilon^* \langle \upsilon g_{\text{right}} \rangle \quad \text{in } L^2(\Omega_{\text{slice}}),$$

where we assume that the eigenfunctions $\upsilon$ and $\upsilon^*$ are normalized such that $\langle \upsilon \upsilon^* \rangle = 1$. Note that the scalar product $\langle \upsilon g_{\text{right}} \rangle$ of positive functions does not vanish. We get

$$(3.24) \qquad \begin{aligned} \lim_{n \to \infty} \int f(\omega_\Lambda) \, \mathbb{P}^{(n)}(d\omega) &= \frac{\langle g^f_{[0,l]} \upsilon^* \rangle \langle \upsilon g_{\text{right}} \rangle}{\langle g^1_{[0,l]} \upsilon^* \rangle \langle \upsilon g_{\text{right}} \rangle} \\ &= \frac{\langle g^f_{[0,l]} \upsilon^* \rangle}{\langle g^1_{[0,l]} \upsilon^* \rangle} = \int f(\omega_\Lambda) \, \mathbb{P}(d\omega), \end{aligned}$$

using the definition (2.41) of $\mathbb{P}$ in the last step. This proves the claim (3.20). To prove the statement for $\tilde{\mathbb{Q}}^{(n)} \to \tilde{\mathbb{Q}}$, consider a bounded measurable function $g((\tilde{x}_e)_{e \in F})$, depending on finitely many components $(\tilde{x}_e)_{e \in F}$, $F \subset \overline{E}$ being finite. Consequently, using the second part of Lemma 3.3, one has for all large $n$:

$$
\begin{aligned}
\int g(\tilde{x}_e^{(n)}; e \in F)\, \tilde{\mathbb{Q}}^{(n)}(d\tilde{x}) &= \int g(\tilde{x}_e^{(n)}(\omega^{(n)}); e \in F)\, \mathbb{P}^{(n)}(d\omega) \\
&\stackrel{n \to \infty}{\longrightarrow} \int g(\tilde{x}_e(\omega); e \in F)\, \mathbb{P}(d\omega) \\
&= \int g(\tilde{x}_e; e \in F)\, \tilde{\mathbb{Q}}(d\tilde{x}),
\end{aligned}
\tag{3.25}
$$

which completes the proof of the lemma. $\square$

LEMMA 3.5 (Exponential decay of the edge weights). *There exist positive constants $c_4, c_5, c_6$ depending only on $\overline{G}$ and $a$ such that for all $n \in \mathbb{N}$ and all edges $e$ of $\overline{G}^{(n)}$, we have*

$$
\mathbb{P}^{(n)}(\tilde{x}_e^{(n)} > e^{-c_4|e|}) \le c_5 e^{-c_6|e|}, \tag{3.26}
$$

*uniformly in $n$. On the infinite ladder $\overline{G}$, we have the similar bound*

$$
\mathbb{P}(\tilde{x}_e > e^{-c_4|e|}) \le c_5 e^{-c_6|e|}. \tag{3.27}
$$

PROOF. We denote by $P_{\mathbf{0}}^{(n)}$ the distribution of the edge-reinforced random walk on the finite graph $\overline{G}^{(n)}$. For $e \in \overline{E}$ and $t \in \mathbb{N}$, let $k_t(e)$ denote the number of times the reinforced random walker traverses the edge $e$ up to time $t$. By Theorem 1.2 of [10] (see also Theorem 1.2 of [7]), there exist positive constants $c_4, c_5, c_6$ depending only on $G$ and $a$ such that for all $n \in \mathbb{N}$ and all edges $e \in \overline{E}^{(n)}$, we have

$$
P_{\mathbf{0}}^{(n)}\left(\lim_{t \to \infty} \frac{k_t(e)}{k_t(e_0^*)} > e^{-c_4|e|}\right) \le c_5 e^{-c_6|e|}, \tag{3.28}
$$

uniformly in $n$. By Theorem 1 of [6], the limit $\lim_{t \to \infty} k_t(e)/t$ exists $P_{\mathbf{0}}^{(n)}$-a.s. and is strictly positive; the distribution of the limiting vector $\lim_{t \to \infty} (k_t(e)/t)_{e \in \overline{E}^{(n)}}$ equals $\mathbb{Q}^{(n)}$, the distribution of the random environment on $\overline{G}^{(n)}$. Furthermore, by (3.19), $x_e^{(n)}/x_{e_0^*}^{(n)} = \tilde{x}_e^{(n)}$. Hence, the left-hand side of (3.28) equals

$$
\begin{aligned}
P_{\mathbf{0}}^{(n)}\left(\lim_{t \to \infty} \frac{k_t(e)/t}{k_t(e_0^*)/t} > e^{-c_4|e|}\right) &= \mathbb{Q}^{(n)}\left(\frac{x_e^{(n)}}{x_{e_0^*}^{(n)}} > e^{-c_4|e|}\right) \\
&= \mathbb{P}^{(n)}(\tilde{x}_e^{(n)}(\omega^{(n)}) > e^{-c_4|e|}).
\end{aligned}
\tag{3.29}
$$

We used Lemma 3.3 in the last step. Thus, the bound (3.28) and equation (3.29) imply the estimate (3.26).

Taking the limit as $n \to \infty$ in (3.26), the claim (3.27) follows from Lemma 3.4. □

LEMMA 3.6 (Normalizability). *The infinite series $\sum_{e \in \overline{E}} \tilde{x}_e$ is $\mathbb{P}$-almost surely finite.*

PROOF. This follows from (3.27) by a Borel–Cantelli argument. □

PROOF OF THEOREMS 2.4 AND 2.21. The sum $\sum_{e \in \overline{E}} \tilde{x}_e$ is $\mathbb{P}$-almost surely finite by Lemma 3.6. Hence, $x(\omega)$, as defined in (2.48), is $\mathbb{P}$-almost surely well defined, and clearly, $\sum_{e \in \overline{E}} x_e(\omega) = 1$.

Let $\pi = (v_0 = \mathbf{0}, v_1, \ldots, v_k)$ be a path in $\overline{G}$. For $n > k$, the random walker cannot visit the ends of the finite graph $\overline{G}^{(n)}$ up to time $k$. Hence, the probability to follow the path $\pi$ up to time $k$ agrees for the reinforced random walker on $\overline{G}$ and on $\overline{G}^{(n)}$:

$$(3.30) \qquad P_{\mathbf{0}}((X_s)_{s=0,\ldots,k} = \pi) = P_{\mathbf{0}}^{(n)}((X_s)_{s=0,\ldots,k} = \pi).$$

Recall that the edge-reinforced random walk on $\overline{G}^{(n)}$ can be represented as a mixture of reversible Markov chains with mixing measure $\mathbb{Q}^{(n)}$. Hence, using Lemmas 3.3 and 3.4,

$$P_{\mathbf{0}}^{(n)}((X_s)_{s=0,\ldots,k} = \pi) = \int Q_{\mathbf{0},x^{(n)}}((X_s)_{s=0,\ldots,k} = \pi) \mathbb{Q}^{(n)}(dx^{(n)})$$

$$= \int Q_{\mathbf{0},\tilde{x}^{(n)}(\omega^{(n)})}((X_s)_{s=0,\ldots,k} = \pi) \mathbb{P}^{(n)}(d\omega^{(n)})$$

$$(3.31) \qquad \overset{n \to \infty}{\longrightarrow} \int Q_{\mathbf{0},\tilde{x}(\omega)}((X_s)_{s=0,\ldots,k} = \pi) \mathbb{P}(d\omega)$$

$$(3.32) \qquad = \int Q_{\mathbf{0},\tilde{x}}((X_s)_{s=0,\ldots,k} = \pi) \tilde{\mathbb{Q}}(d\tilde{x});$$

recall that $\tilde{\mathbb{Q}}$ is the law of $\tilde{x}$ with respect to $\mathbb{P}$.

Now, any random walk distribution is uniquely determined by its values on the events $A = \{(X_s)_{s=0,\ldots,k} = \pi\}$. Thus, (3.30) and (3.32) imply that edge-reinforced random walk on $\overline{G}$ has the same distribution as a random walk in a random environment with weights having the distribution $\tilde{\mathbb{Q}}$. By (2.46), the environment is normalized so that $\tilde{x}_{e_0^*} = 1$.

Moreover, the weights $\tilde{x}_e$ and $x_e$ are proportional by (3.19); hence the Markov chain probabilities $Q_{\mathbf{0},\tilde{x}}$ and $Q_{\mathbf{0},x}$ coincide. Thus, it follows from (3.30) and (3.31) that edge-reinforced random walk on $\overline{G}$ has also the same distribution as a random walk in a random environment $\mathbb{Q}$, given by weights $x(\omega)$, where $\omega$ has the distribution $\mathbb{P}$.

By Theorem (7) of [3], the mixing measure $\mathbb{Q}$ is uniquely determined, because the reinforced random walk on $\overline{G}$ is recurrent as was shown in Lemma 5.2 of [7] and Proposition 5.1 of [10]. We conclude that (2.8) holds. □

PROOF OF THEOREM 2.5. By the conversion (3.19), $x_e^{(n)}/x_{e_0^*}^{(n)} = \tilde{x}_e^{(n)}$ and $x_e/x_{e_0^*} = \tilde{x}_e$. Using Lemma 3.3 and Theorem 2.21, the claims follow immediately from Lemma 3.5. □

Let us summarize the key statements in Lemma 3.3 and Theorem 2.21.
*Connection between the random environment and the Gibbs measures*:
$\mathbb{Q}$ equals the law of $x(\omega)$ if $\omega$ has law $\mathbb{P}$.
$\mathbb{Q}^{(n)}$ equals the law of $x^{(n)}(\omega^{(n)})$ if $\omega^{(n)}$ has law $\mathbb{P}^{(n)}$.

3.2. *Properties of the random environment.*

PROOF OF THEOREM 2.7. First, we show that for any finite path $\rho = (v_k, \ldots, v_l)$ with $l \geq k$, we have

$$
\begin{aligned}
(3.33) \quad & P((X_s)_{s=k,\ldots,l} = \rho | (X_s)_{s=0,\ldots,k} = \pi) \\
& = \int_\Delta Q_{v_k,x}((X_s)_{s=k,\ldots,l} = \rho) \frac{Q_{\mathbf{0},x}((X_s)_{s=0,\ldots,k} = \pi)}{P_{\mathbf{0}}((X_s)_{s=0,\ldots,k} = \pi)} \mathbb{Q}(dx).
\end{aligned}
$$

Let $\pi\rho$ denote the concatenation of $\pi$ and $\rho$. By the Markov property for $Q_{\cdot,x}$ and the representation of the edge-reinforced random walk as a mixture of Markov chains [identity (2.8) in Theorem 2.4], the right-hand side of (3.33) equals

$$
\begin{aligned}
& \frac{\int_\Delta Q_{\mathbf{0},x}((X_s)_{s=0,\ldots,l} = \pi\rho) \mathbb{Q}(dx)}{P_{\mathbf{0}}((X_s)_{s=0,\ldots,k} = \pi)} \\
(3.34) \quad & = \frac{P_{\mathbf{0}}((X_s)_{s=0,\ldots,l} = \pi\rho)}{P_{\mathbf{0}}((X_s)_{s=0,\ldots,k} = \pi)} \\
& = P_{\mathbf{0}}((X_s)_{s=k,\ldots,l} = \rho | (X_s)_{s=0,\ldots,k} = \pi).
\end{aligned}
$$

Thus, (3.33) holds.

Since the distribution of $(X_{k+t})_{t \in \mathbb{N}_0}$ is uniquely determined by its values on events of the form $\{(X_s)_{s=k\ldots l} = \rho\}$, (3.33) generalizes to

$$
\begin{aligned}
(3.35) \quad & P((X_{k+t})_{t \in \mathbb{N}_0} \in A | (X_s)_{s=0,\ldots,k} = \pi) \\
& = \int_\Delta Q_{v_k,x}((X_{k+t})_{t \in \mathbb{N}_0} \in A) \frac{Q_{\mathbf{0},x}((X_s)_{s=0,\ldots,k} = \pi)}{P_{\mathbf{0}}((X_s)_{s=0,\ldots,k} = \pi)} \mathbb{Q}(dx).
\end{aligned}
$$

This proves the claim (2.13) for the measure $\mathbb{Q}_\pi$ with the Radon–Nikodym derivative (2.14). In particular, taking $A = \overline{V}^{\mathbb{N}_0}$ in (2.13), we see that $\mathbb{Q}_\pi$ is a probability measure. The uniqueness of $\mathbb{Q}_\pi$ follows immediately from the uniqueness of $\mathbb{Q}$ stated in Theorem 2.4. $\square$

PROOF OF THEOREM 2.8. By Lemma 3.4, the finite-dimensional marginals of $\tilde{\mathbb{Q}}^{(n)}$ converge weakly to the corresponding marginals of $\tilde{\mathbb{Q}}$. By Theorem 10.1 in Chapter 3 of [11], there exists a coupling $((\hat{x}_e^{(n)})_{e \in \overline{E}}, n \in \mathbb{N}, (\hat{x}_e)_{e \in \overline{E}})$ with a coupling measure $\hat{\mathbb{Q}}$, such that

$$(3.36) \qquad \lim_{n \to \infty} \hat{x}_e^{(n)} = \hat{x}_e \qquad \text{pointwise for all } e \in \overline{E},$$

$(\hat{x}_e)_{e \in \overline{E}}$ has the law $\tilde{\mathbb{Q}}$, and $(\hat{x}_e^{(n)})_{e \in \overline{E}}$ has the law $\tilde{\mathbb{Q}}^{(n)}$. More precisely, we set $\hat{x}_e^{(n)} = 0$ whenever $e$ is not an edge in $\overline{G}^{(n)}$, and we let $(\hat{x}_e^{(n)})_{e \in \overline{E}^{(n)}}$ have the law $\tilde{\mathbb{Q}}^{(n)}$.

We claim that

$$(3.37) \qquad \lim_{n \to \infty} \sum_{e \in \overline{E}} \hat{x}_e^{(n)} = \sum_{e \in \overline{E}} \hat{x}_e, \qquad \hat{\mathbb{Q}}\text{-a.s.}$$

This can be seen as follows. Recall that $\tilde{x}_e = x_e/x_{e_0^*}$. Hence, from (2.10) of Theorem 2.5, we know that for every edge $e$,

$$(3.38) \qquad \hat{\mathbb{Q}}(\hat{x}_e > e^{-c_4|e|}) \le c_5 e^{-c_6|e|}.$$

Let $\varepsilon > 0$. Take $m$ so large that $\sum_{e \in \overline{E} \setminus \overline{E}^{(m)}} c_5 e^{-c_6|e|} < \varepsilon/2$ and $\sum_{e \in \overline{E} \setminus \overline{E}^{(m)}} e^{-c_4|e|} < \varepsilon/2$. We estimate

$$(3.39) \qquad \begin{aligned} & \left| \sum_{e \in \overline{E}} \hat{x}_e^{(n)} - \sum_{e \in \overline{E}} \hat{x}_e \right| \\ & \le \sum_{e \in \overline{E}^{(m)}} |\hat{x}_e^{(n)} - \hat{x}_e| + \sum_{e \in \overline{E} \setminus \overline{E}^{(m)}} [|\hat{x}_e^{(n)}| + |\hat{x}_e|]. \end{aligned}$$

The convergence (3.36) implies that the first sum on the right-hand side of (3.39) converges to 0 $\hat{\mathbb{Q}}$-a.s. as $n \to \infty$. Therefore, we conclude

$$(3.40) \qquad \begin{aligned} & \hat{\mathbb{Q}}\left( \limsup_{n \to \infty} \left| \sum_{e \in \overline{E}} \hat{x}_e^{(n)} - \sum_{e \in \overline{E}} \hat{x}_e \right| > \varepsilon \right) \\ & \le \hat{\mathbb{Q}}\left( \limsup_{n \to \infty} \sum_{e \in \overline{E} \setminus \overline{E}^{(m)}} [|\hat{x}_e^{(n)}| + |\hat{x}_e|] > \varepsilon \right) \\ & \le \sum_{e \in \overline{E} \setminus \overline{E}^{(m)}} [\hat{\mathbb{Q}}(\hat{x}_e^{(n)} > e^{-c_4|e|}) + \hat{\mathbb{Q}}(\hat{x}_e > e^{-c_4|e|})] < \varepsilon. \end{aligned}$$

This implies the claim (3.37).

As a consequence of our normalization $\hat{x}_{e_0} = 1$, we know $\sum_{e \in \overline{E}} \hat{x}_e \geq 1$. Hence, (3.36) and (3.37) imply that we have for all $e \in \overline{E}$

$$(3.41) \quad \lim_{n \to \infty} \frac{\hat{x}_e^{(n)}}{\sum_{e' \in \overline{E}^{(n)}} \hat{x}_{e'}^{(n)}} = \lim_{n \to \infty} \frac{\hat{x}_e^{(n)}}{\sum_{e' \in \overline{E}} \hat{x}_{e'}^{(n)}} = \frac{\hat{x}_e}{\sum_{e' \in \overline{E}} \hat{x}_{e'}}, \quad \hat{\mathbb{Q}}\text{-a.s.}$$

We know that the $\hat{\mathbb{Q}}$-distribution of $(\hat{x}_e^{(n)} / \sum_{e' \in \overline{E}^{(n)}} \hat{x}_{e'}^{(n)})_{e \in \overline{E}^{(n)}}$ equals $\mathbb{Q}^{(n)}$ by Lemma 3.3, and the $\hat{\mathbb{Q}}$-distribution of $(\hat{x}_e / \sum_{e' \in \overline{E}} \hat{x}_{e'})_{e \in \overline{E}}$ equals $\mathbb{Q}$ by Theorem 2.21. Hence, (3.41) implies that the finite-dimensional marginals of $\mathbb{Q}^{(n)}$ converge weakly to the corresponding marginals of $\mathbb{Q}$. □

PROOF OF THEOREM 2.9. Since the sum of finitely many random variables with exponential tails has again exponential tails, Theorem 2.2 of [10] (see also Theorem 2.3 of [7]) implies the following tail estimate in finite volume: There exist positive constants $c_7(a), c_8(a)$, depending only on $a$, such that one has for all $n \in \mathbb{N}$, $i, j \in \{0, 1, \ldots, n\}$ with $|i - j| \leq 1$ and $M > 0$,

$$(3.42) \quad \mathbb{Q}^{(n)}\left[\left|\ln \frac{x_e}{x_{e'}}\right| \geq M\right] \leq c_7 e^{-c_8 M},$$

whenever $e, e'$ are edges on level $i, j$, respectively.

Now, let $e, e'$ be edges on arbitrary levels $i \neq j$, respectively. Then, we can write $\ln(x_e/x_{e'})$ as a sum of $|i - j|$ terms of the form $\ln(x_f/x_{f'})$ with $f$ and $f'$ edges on neighboring levels. Hence,

$$(3.43) \quad \mathbb{Q}^{(n)}\left[\left|\ln \frac{x_e}{x_{e'}}\right| \geq M\right] \leq c_7 |i - j| e^{-c_8 M / |i-j|}$$

follows. Note that the constants $c_7$ and $c_8$ are independent of $n$. By Theorem 2.8, we can take the limit as $n \to \infty$ in the inequality (3.43); note that the distributions of all $\ln(x_e/x_{e'})$, $e \neq e'$, with respect to the limit law $\mathbb{Q}$ are continuous. This yields the claim (2.16) and completes the proof of the theorem. □

PROOF OF COROLLARY 2.6. Fix $f \in \overline{E}$. For $e \in \overline{E}$, we define the event

$$(3.44) \quad A_e := \{x_e > e^{-c_4|e|/2} x_f\}$$

with $c_4$ as in Theorem 2.5. By (2.46) the weights $\tilde{x}_e$ are normalized so that $\tilde{x}_{e_0^*} = 1$. Furthermore, $\tilde{x}_e$ and $x_e$ are proportional with $\tilde{x}_e = x_e/x_{e_0^*}$. An application of the estimate (2.10) from Theorem 2.5 and the bound (2.16) from Theorem 2.9 yield

$$(3.45) \quad \begin{aligned} \mathbb{Q}(A_e) &\leq \mathbb{Q}(x_e > e^{-c_4|e|} x_{e_0^*}) + \mathbb{Q}(x_{e_0^*} > e^{c_4|e|/2} x_f) \\ &\leq c_5 e^{-c_6|e|} + c_9 e^{-c_{10}|e|} \end{aligned}$$

with constants $c_9(f) > 0$ and $c_{10}(f) > 0$. Thus, $\sum_{e \in \overline{E}} \mathbb{Q}(A_e) < \infty$, and the Borel–Cantelli lemma implies that $\mathbb{Q}$-a.s., the event $A_e$ holds for at most finitely many $e \in \overline{E}$, that is, $\mathbb{Q}$-a.s., there exists $n_0$ such that for all $n \geq n_0$, the claim (2.11) holds. $\square$

3.3. *Asymptotic properties of the reinforced random walk.*

PROOF OF THEOREM 2.1. Let $v \in \overline{V}$ be a vertex and $x \in \Delta$. Since $\sum_{e \in \overline{E}} x_e = 1$, we know $\sum_{v \in \overline{V}} x_v = 2$; recall the notation (2.7). Furthermore, with respect to the invariant distribution $\pi = (x_v/2)_{v \in \overline{V}}$, the random walk $Q_{\pi,x}$ in the fixed environment $x$ is reversible; see, for example, Example 4.5 on pages 298–299 of [5]. We conclude that

$$(3.46) \qquad x_{\mathbf{0}} Q_{\mathbf{0},x}(X_t = v) = x_v Q_{v,x}(X_t = \mathbf{0})$$

and therefore

$$(3.47) \qquad Q_{\mathbf{0},x}(X_t = v) = \frac{x_v}{x_{\mathbf{0}}} Q_{v,x}(X_t = \mathbf{0}) \leq \frac{x_v}{x_{\mathbf{0}}} \leq \frac{x_v}{x_{e_0^*}}.$$

Taking now a random environment $x$, we consider the event

$$(3.48) \qquad A_v = \left\{ \frac{x_v}{x_{e_0^*}} \geq k e^{-c_4(|v|-1)} \right\},$$

where $k$ denotes the coordination number of $\overline{G}$, i.e. the maximal number of immediate neighbors that a vertex in $\overline{G}$ can have. We estimate, using the bound (3.27) from Lemma 3.5 and the fact $x_v/x_{e_0^*} = \sum_{e \in \overline{E}: v \in e} \tilde{x}_e$,

$$(3.49) \qquad \begin{aligned} \mathbb{P}(A_v) &\leq \sum_{e \in \overline{E}: v \in e} \mathbb{P}(\tilde{x}_e \geq e^{-c_4(|v|-1)}) \\ &\leq \sum_{e \in \overline{E}: v \in e} \mathbb{P}(\tilde{x}_e \geq e^{-c_4|e|}) \\ &\leq k c_5 e^{-c_6(|v|-1)}; \end{aligned}$$

note that $|e| \geq |v| - 1$ if $v \in e$. Using the representation (2.8) of $P_{\mathbf{0}}$ as a mixture of the $Q_{\mathbf{0},x}$, we get the bound

$$(3.50) \qquad \begin{aligned} P_{\mathbf{0}}(X_t = v) &= \int_\Omega Q_{\mathbf{0},x(\omega)}(X_t = v) \, \mathbb{P}(d\omega) \\ &= \int_{A_v^c} Q_{\mathbf{0},x}(X_t = v) \, d\mathbb{P} + \int_{A_v} Q_{\mathbf{0},x}(X_t = v) \, d\mathbb{P} \\ &\leq \int_{A_v^c} \frac{x_v}{x_{e_0^*}} \, d\mathbb{P} + \mathbb{P}(A_v) \\ &\leq k e^{-c_4(|v|-1)} + k c_5 e^{-c_6(|v|-1)}. \end{aligned}$$

The claim (2.1) of Theorem 2.1 follows from this bound by summation over all vertices $v$ with $|v| \geq n$. $\square$

PROOF OF COROLLARY 2.2. By Theorem 2.1, we know that for all $c_3 > 0$ and $t \geq 2$, the following holds:

$$(3.51) \quad P_{\mathbf{0}}\left(\max_{s=0,\ldots,t} |X_s| > c_3 \ln t\right) \leq \sum_{s=1}^{t} P_{\mathbf{0}}(|X_s| > c_3 \ln t) \leq t c_1 e^{-c_2 c_3 \ln t}.$$

We choose $c_3$ large enough that $1 - c_2 c_3 \leq -2$. Then

$$(3.52) \quad P_{\mathbf{0}}\left(\max_{s=0,\ldots,t} |X_s| > c_3 \ln t\right) \leq c_1 t^{-2};$$

in particular the probabilities in (3.52) are summable over $t$. Hence, by the Borel-Cantelli lemma, we know that $P_{\mathbf{0}}$-a.s. the claim (2.2) holds. $\square$

PROOF OF THEOREM 2.12. Let $x \in \Delta$. Recall that $(x_v/2)_{v \in \overline{V}}$ is a reversible probability distribution for the Markov chain $Q_{\cdot,x}$. All closed paths in $\overline{G} = \mathbb{N}_0 \times G$ have an even length because $G$ is acyclic. Hence, all states $v \in \overline{V}$ have period 2. Furthermore, the Markov chain $Q_{\cdot,x}$ is irreducible. Consequently, by the convergence theorem (5.7) on page 315 of [5], (2.23) and (2.24) follow.

To prove (2.25), let $A \subseteq \overline{V}^{\mathbb{N}_0}$ be measurable. Note that the convergence in (2.23) holds also in the $l^1$-norm on $\mathbb{R}^{\overline{V}}$ by the discrete version of Scheffé's theorem. Hence, we can exchange infinite sum and limit in the following calculation:

$$\lim_{s \to \infty} Q_{\mathbf{0},x}((X_{2s+t})_{t \in \mathbb{N}_0} \in A)$$

$$= \lim_{s \to \infty} \sum_{v \in \overline{V}} Q_{\mathbf{0},x}((X_{2s+t})_{t \in \mathbb{N}_0} \in A | X_{2s} = v) Q_{\mathbf{0},x}(X_{2s} = v)$$

$$(3.53) \quad = \lim_{s \to \infty} \sum_{v \in \overline{V}} Q_{v,x}((X_t)_{t \in \mathbb{N}_0} \in A) Q_{\mathbf{0},x}(X_{2s} = v)$$

$$= \sum_{v \in \overline{V}} Q_{v,x}((X_t)_{t \in \mathbb{N}_0} \in A) \lim_{s \to \infty} Q_{\mathbf{0},x}(X_{2s} = v)$$

$$= Q_{x_{\text{even}},x}(A).$$

By the same argument, the last claim (2.26) follows. $\square$

PROOF OF THEOREMS 2.3 AND 2.10. We combine the representation of the reinforced random walk as a mixture of Markov chains (Theorem 2.4)

and the convergence to equilibrium, conditioned on the environment (Theorem 2.12). By the dominated convergence theorem, this yields

$$\lim_{t\to\infty} P_{\mathbf{0}}(X_{2t} = v) = \int_\Delta \lim_{t\to\infty} Q_{\mathbf{0},x}(X_{2t} = v)\,\mathbb{Q}(dx)$$
$$(3.54) \qquad = \int_\Delta x_{\text{even}}(v)\,\mathbb{Q}(dx) =: \mu_{\text{even}}(v)$$

for all $v \in \overline{V}$. Since $x_{\text{even}}$ is a probability function, by the monotone convergence theorem, $\mu_{\text{even}}$ is a probability function as well.

In order to prove the upper bound for $\mu_{\text{even}}(v)$, we use (3.54) and Theorem 2.1:

$$(3.55)\quad \mu_{\text{even}}(v) = \lim_{t\to\infty} P_{\mathbf{0}}(X_{2t} = v) \le \limsup_{t\to\infty} P_{\mathbf{0}}(|X_{2t}| \ge |v|) \le c_1 e^{-c_2|v|}.$$

The convergence of the distribution of $X_{2t+1}$ and the claims for $\mu_{\text{odd}}$ are proved analogously. □

PROOF OF THEOREM 2.11. The claim follows from Theorems 2.4 and 2.12 and the dominated convergence theorem. □

**Acknowledgments.** It is a great pleasure to thank Lincoln Chayes for his inspiring lectures, which influenced the style of this paper. We are grateful to the University of Leiden and the University of California at Los Angeles for their hospitality.

MATHEMATICAL INSTITUTE  
UNIVERSITY OF MUNICH  
THERESIENSTR. 39  
D-80333 MUNICH  
GERMANY  
E-MAIL: merkl@mathematik.uni-muenchen.de  
URL: http://www.mathematik.uni-muenchen.de/~merkl  

ZENTRUM MATHEMATIK  
BEREICH M5  
TECHNISCHE UNIVERSITÄT MÜNCHEN  
D-85747 GARCHING BEI MÜNCHEN  
GERMANY  
E-MAIL: srolles@ma.tum.de  
URL: http://www-m5.ma.tum.de/pers/srolles/